# Rate of convergence of predictive distributions for dependent data

PATRIZIA BERTI[1], IRENE CRIMALDI[2], LUCA PRATELLI[3] and
PIETRO RIGO[4]

[1]*Dipartimento di Matematica Pura ed Applicata "G. Vitali", Universita' di Modena e Reggio-Emilia, via Campi 213/B, 41100 Modena, Italy. E-mail: patrizia.berti@unimore.it*
[2]*Dipartimento di Matematica, Universita' di Bologna, Piazza di Porta San Donato 5, 40126 Bologna, Italy. E-mail: crimaldi@dm.unibo.it*
[3]*Accademia Navale, viale Italia 72, 57100 Livorno, Italy. E-mail: pratel@mail.dm.unipi.it*
[4]*Dipartimento di Economia Politica e Metodi Quantitativi, Universita' di Pavia, via S. Felice 5, 27100 Pavia, Italy. E-mail: prigo@eco.unipv.it*

This paper deals with empirical processes of the type

$$C_n(B) = \sqrt{n}\{\mu_n(B) - P(X_{n+1} \in B \mid X_1, \ldots, X_n)\},$$

where $(X_n)$ is a sequence of random variables and $\mu_n = (1/n)\sum_{i=1}^n \delta_{X_i}$ the empirical measure. Conditions for $\sup_B |C_n(B)|$ to converge stably (in particular, in distribution) are given, where $B$ ranges over a suitable class of measurable sets. These conditions apply when $(X_n)$ is exchangeable or, more generally, conditionally identically distributed (in the sense of Berti *et al.* [*Ann. Probab.* **32** (2004) 2029–2052]). By such conditions, in some relevant situations, one obtains that $\sup_B |C_n(B)| \xrightarrow{P} 0$ or even that $\sqrt{n}\sup_B |C_n(B)|$ converges a.s. Results of this type are useful in Bayesian statistics.

## 1. Introduction and motivations

A number of real problems reduce to the evaluation of the *predictive distribution*

$$a_n(\cdot) = P(X_{n+1} \in \cdot | X_1, \ldots, X_n)$$

for a sequence $X_1, X_2, \ldots$ of random variables. Here, we focus on those situations where $a_n$ cannot be calculated in closed form and one decides to estimate it based on the available data $X_1, \ldots, X_n$. Related references are [1–3, 5, 6, 8, 10, 15, 18, 20].

For notational reasons, it is convenient to work in coordinate probability space. Accordingly, we fix a measurable space $(S, \mathcal{B})$ and a probability $P$ on $(S^\infty, \mathcal{B}^\infty)$, and we let $X_n$ be the $n$th canonical projection on $(S^\infty, \mathcal{B}^\infty, P)$, $n \geq 1$. We also let

$$\mathcal{G}_n = \sigma(X_1, \ldots, X_n) \quad \text{and} \quad X = (X_1, X_2, \ldots).$$

Since we are concerned with predictive distributions, it is reasonable to make some (qualitative) assumptions about them. In [6], $X$ is said to be *conditionally identically distributed* (c.i.d.) when

$$E(I_B(X_k)|\mathcal{G}_n) = E(I_B(X_{n+1})|\mathcal{G}_n) \quad \text{a.s. for all } B \in \mathcal{B} \text{ and } k > n \geq 0,$$

where $\mathcal{G}_0$ is the trivial $\sigma$-field. Thus, at each time $n \geq 0$, the future observations $(X_k : k > n)$ are identically distributed given the past $\mathcal{G}_n$. In a sense, this is a weak form of exchangeability. In fact, $X$ is exchangeable if and only if it is stationary and c.i.d., and various examples of non-exchangeable c.i.d. sequences are available.

In the sequel, $X = (X_1, X_2, \ldots)$ is a c.i.d. sequence of random variables.

In that case, a sound estimate of $a_n$ is the *empirical distribution*

$$\mu_n = \frac{1}{n} \sum_{i=1}^n \delta_{X_i}.$$

The choice of $\mu_n$ can be defended as follows. Let $\mathcal{D} \subset \mathcal{B}$ and let $\|\cdot\|$ denote the sup-norm on $\mathcal{D}$. Suppose also that $\mathcal{D}$ is countably determined, as defined in Section 2. (The latter is a mild condition, only needed to handle measurability issues.) Then

$$\|\mu_n - a_n\| = \sup_{B \in \mathcal{D}} |\mu_n(B) - a_n(B)| \xrightarrow{\text{a.s.}} 0, \tag{1}$$

provided ($X$ is c.i.d. and) $\mu_n$ converges uniformly on $\mathcal{D}$ with probability 1; see [5]. For instance, $\|\mu_n - a_n\| \xrightarrow{\text{a.s.}} 0$ whenever $X$ is exchangeable and $\mathcal{D}$ is a Glivenko–Cantelli class. Also, $\|\mu_n - a_n\| \xrightarrow{\text{a.s.}} 0$ if $S = \mathbb{R}$, $\mathcal{D} = \{(-\infty, t] : t \in \mathbb{R}\}$, and $X_1$ has a discrete distribution or $\inf_{\varepsilon > 0} \liminf_n P(|X_{n+1} - X_n| < \varepsilon) = 0$; see [4].

To sum up, under mild assumptions, $\mu_n$ is a consistent estimate of $a_n$ (with respect to uniform distance) for c.i.d. data. This is in line with de Finetti [10] in the particular case of exchangeable indicators.

Taking (1) as a starting point, the next step is to investigate the convergence rate, that is, to investigate whether $\alpha_n \|\mu_n - a_n\|$ converges in distribution, possibly to a null limit, for suitable constants $\alpha_n > 0$. This is precisely the purpose of this paper.

A first piece of information on the convergence rate of $\|\mu_n - a_n\|$ can be obtained as follows. For $B \in \mathcal{B}$, define

$$\mu(B) = \limsup_n \mu_n(B),$$

$$W_n(B) = \sqrt{n} \{\mu_n(B) - \mu(B)\}.$$

By the SLLN for c.i.d. sequences, $\mu_n(B) \xrightarrow{a.s.} \mu(B)$; see [6]. Hence, for fixed $n \geq 0$ and $B \in \mathcal{B}$, one obtains

$$E(\mu(B)|\mathcal{G}_n) = \lim_k E(\mu_k(B)|\mathcal{G}_n) = \lim_k \frac{1}{k} \sum_{i=n+1}^k E(I_B(X_i)|\mathcal{G}_n)$$
$$= E(I_B(X_{n+1})|\mathcal{G}_n) = a_n(B) \quad \text{a.s.}$$

In turn, this implies that $\sqrt{n}\{\mu_n(B) - a_n(B)\} = E(W_n(B)|\mathcal{G}_n)$ a.s., so

$$\|\mu_n - a_n\| \leq \frac{1}{\sqrt{n}} \sup_{B \in \mathcal{D}} E(|W_n(B)||\mathcal{G}_n) \leq \frac{1}{\sqrt{n}} E(\|W_n\||\mathcal{G}_n) \quad \text{a.s.}$$

If $\sup_n E\|W_n\|^k < \infty$ for some $k \geq 1$, it then follows that

$$E\{(\alpha_n\|\mu_n - a_n\|)^k\} \leq \left(\frac{\alpha_n}{\sqrt{n}}\right)^k E\|W_n\|^k \to 0 \quad \text{whenever } \frac{\alpha_n}{\sqrt{n}} \to 0.$$

Even if obvious, this fact is potentially useful since

$$\sup_n E\|W_n\|^k < \infty \quad \text{for all } k \geq 1, \text{ if } X \text{ is exchangeable}, \tag{2}$$

for various choices of $\mathcal{D}$; see Remark 3. In particular, (2) holds if $\mathcal{D}$ is finite.

The intriguing case, however, is $\alpha_n = \sqrt{n}$. For each $B \in \mathcal{B}$ and probability $Q$ on $(S^\infty, \mathcal{B}^\infty)$, write

$$C_n^Q(B) = E_Q(W_n(B)|\mathcal{G}_n) \quad \text{and}$$
$$C_n(B) = C_n^P(B) = \sqrt{n}\{\mu_n(B) - a_n(B)\}.$$

In Theorem 3.3 of [6], the asymptotic behavior of $C_n(B)$ is investigated for *fixed* $B$. Here, instead, we are interested in

$$\|C_n\| = \sup_{B \in \mathcal{D}} |C_n(B)| = \sqrt{n}\|\mu_n - a_n\|.$$

Our main result (Theorem 1) is the following. Fix a random probability measure $N$ on $\mathbb{R}$ and a probability $Q$ on $(S^\infty, \mathcal{B}^\infty)$ such that

$$\|C_n^Q\| \to N \quad \text{stably under } Q \quad \text{and}$$
$$\|W_n\| \text{ is uniformly integrable under both } P \text{ and } Q.$$

Then,

$$\|C_n\| \to N \quad \text{stably whenever } P \ll Q. \tag{3}$$

A remarkable particular case is $N = \delta_0$. Suppose, in fact, that for some $Q$, one has $\|C_n^Q\| \xrightarrow{Q} 0$ and $\|W_n\|$ uniformly integrable under $P$ and $Q$. Then,

$$\|C_n\| \xrightarrow{P} 0 \qquad \text{whenever } P \ll Q.$$

Stable convergence (in the sense of Rényi) is a stronger form of convergence in distribution. The definition is recalled in Section 2.

In general, one cannot dispense with the uniform integrability condition. However, this condition is often true. For instance, $\|W_n\|$ is uniformly integrable (under $P$ and $Q$) provided $\mathcal{D}$ meets (2) and $X$ is exchangeable (under $P$ and $Q$).

To make (3) concrete, a large list of reference probabilities $Q$ is needed. Various examples are available in the Bayesian nonparametrics framework; see, for example, [16] and references therein. The most popular is perhaps the Ferguson–Dirichlet law, denoted by $Q_0$. If $P = Q_0$, then $X$ is exchangeable and

$$a_n(B) = \frac{\alpha P(X_1 \in B) + n\mu_n(B)}{\alpha + n} \qquad \text{a.s. for some constant } \alpha > 0.$$

Since $\|\mu_n - a_n\| \leq (\alpha/n)$ when $P = Q_0$, something more than $\|C_n\| \xrightarrow{P} 0$ can be expected in the case $P \ll Q_0$. Indeed, we prove that

$$n\|\mu_n - a_n\| = \sqrt{n}\|C_n\| \qquad \text{converges a.s.}$$

whenever $P \ll Q_0$ with a density satisfying a certain condition; see Theorem 2 and Corollary 5.

One more example should be mentioned. Let $X_n = (Y_n, Z_n)$, where $Z_n > 0$ and

$$P(Y_{n+1} \in B | \mathcal{G}_n) = \frac{\alpha P(Y_1 \in B) + \sum_{i=1}^n Z_i I_B(Y_i)}{\alpha + \sum_{i=1}^n Z_i} \qquad \text{a.s.}$$

for some constant $\alpha > 0$. Under some conditions, $X$ is c.i.d. (but not necessarily exchangeable), $\|W_n\|$ is uniformly integrable and $\|C_n\|$ converges stably; see Section 4.

The above material takes a nicer form when the condition $P \ll Q$ can be given a simple characterization. This happens, for instance, if $S = \{x_1, \ldots, x_k, x_{k+1}\}$ is finite, $X$ exchangeable and $P(X_1 = x) > 0$ for all $x \in S$. Then, $P \ll Q_0$ (for some choice of $Q_0$) if and only if

$$(\mu\{x_1\}, \ldots, \mu\{x_k\})$$

has an absolutely continuous distribution with respect to Lebesgue measure. In this particular case, however, a part of our results can also be obtained through the Bernstein–von Mises theorem; see Section 3.

Finally, we make two remarks:

(i) If $X$ is exchangeable, our results apply to Bayesian predictive inference. Suppose, in fact, that $S$ is Polish and $\mathcal{B}$ the Borel $\sigma$-field, so that de Finetti's theorem applies. Then, $P$ is a unique mixture of product probabilities on $\mathcal{B}^\infty$ and the mixing measure is called the *prior distribution* in a Bayesian framework. Now, given $Q$, $P \ll Q$ is just an assumption on the prior distribution. This is plain in the last example where $S = \{x_1, \ldots, x_k, x_{k+1}\}$. In Bayesian terms, such an example can be summarized as follows. For a multinomial statistical model, $\|C_n\| \xrightarrow{P} 0$ if the prior is absolutely continuous with respect to Lebesgue measure, and $\sqrt{n}\|C_n\|$ converges a.s. if the prior density satisfies a certain condition.

(ii) To our knowledge, there is no general representation for the predictive distributions of an exchangeable sequence. Such a representation would be very useful. Even if only partially, results like (3) contribute to filling the gap. As an example, for fixed $B \in \mathcal{B}$, one obtains $a_n(B) = \mu_n(B) + o_P(\frac{1}{\sqrt{n}})$, provided $X$ is exchangeable and $P \ll Q$ for some $Q$ such that $C_n^Q(B) \xrightarrow{Q} 0$ and $W_n(B)$ is uniformly integrable.

## 2. Main results

A few definitions need to be recalled. Let $T$ be a metric space, $\mathcal{B}_T$ the Borel $\sigma$-field on $T$ and $(\Omega, \mathcal{A}, P)$ a probability space. A *random probability measure* on $T$ is a mapping $N$ on $\Omega \times \mathcal{B}_T$ such that: (i) $N(\omega, \cdot)$ is a probability on $\mathcal{B}_T$ for each $\omega \in \Omega$; (ii) $N(\cdot, B)$ is $\mathcal{A}$-measurable for each $B \in \mathcal{B}_T$. Let $(Z_n)$ be a sequence of $T$-valued random variables and $N$ a random probability measure on $T$. Both $(Z_n)$ and $N$ are defined on $(\Omega, \mathcal{A}, P)$. We say that $Z_n$ *converges stably* to $N$ in the case where

$$P(Z_n \in \cdot | H) \to E(N(\cdot)|H) \qquad \text{weakly for all } H \in \mathcal{A} \text{ such that } P(H) > 0.$$

Clearly, if $Z_n \to N$ stably, then $Z_n$ converges in distribution to the probability law $E(N(\cdot))$ (just let $H = \Omega$). Stable convergence has been introduced by Rényi in [17] and subsequently investigated by various authors; see [9] for more information.

Next, we say that $\mathcal{D} \subset \mathcal{B}$ is *countably determined* in the case where, for some fixed countable subclass $\mathcal{D}_0 \subset \mathcal{D}$, one obtains $\sup_{B \in \mathcal{D}_0} |\nu_1(B) - \nu_2(B)| = \sup_{B \in \mathcal{D}} |\nu_1(B) - \nu_2(B)|$ for every pair $\nu_1, \nu_2$ of probabilities on $\mathcal{B}$. A sufficient condition is that for some countable $\mathcal{D}_0 \subset \mathcal{D}$, and for every $\varepsilon > 0$, $B \in \mathcal{D}$ and probability $\nu$ on $\mathcal{B}$, there is $B_0 \in \mathcal{D}_0$ satisfying $\nu(B \Delta B_0) < \varepsilon$. Most classes $\mathcal{D}$ involved in applications are countably determined. For instance, $\mathcal{D} = \{(-\infty, t] : t \in \mathbb{R}^k\}$ and $\mathcal{D} = \{\text{closed balls}\}$ are countably determined if $S = \mathbb{R}^k$ and $\mathcal{B}$ is the Borel $\sigma$-field. As another example, $\mathcal{D} = \mathcal{B}$ is countably determined if $\mathcal{B}$ is countably generated.

We are now in a position to state our main result. Let $N$ be a random probability measure on $\mathbb{R}$, defined on the measurable space $(S^\infty, \mathcal{B}^\infty)$, and let $Q$ be a probability on $(S^\infty, \mathcal{B}^\infty)$.

**Theorem 1.** *Let $\mathcal{D}$ be countably determined. Suppose $\|C_n^Q\| \to N$ stably under $Q$, and $(\|W_n\| : n \geq 1)$ is uniformly integrable under $P$ and $Q$. Then,*

$$\|C_n\| = \sqrt{n}\|\mu_n - a_n\| \to N \qquad \text{stably whenever } P \ll Q.$$

**Proof.** Since $\mathcal{D}$ is countably determined, there are no measurability problems in taking $\sup_{B \in \mathcal{D}}$. In particular, $\|W_n\|$ and $\|C_n\|$ are random variables and $\|C_n\|$ is $\mathcal{G}_n$-measurable. Let $f$ be a version of $\frac{dP}{dQ}$ and $U_n = f - E_Q(f|\mathcal{G}_n)$. Then,

$$C_n(B) = E(W_n(B)|\mathcal{G}_n) = \frac{E_Q(fW_n(B)|\mathcal{G}_n)}{E_Q(f|\mathcal{G}_n)}$$

$$= C_n^Q(B) + \frac{E_Q(U_n W_n(B)|\mathcal{G}_n)}{E_Q(f|\mathcal{G}_n)}, \qquad P\text{-a.s., for each } B \in \mathcal{B}.$$

Letting $M_n = \frac{E_Q(|U_n|\|W_n\| |\mathcal{G}_n)}{E_Q(f|\mathcal{G}_n)}$ and taking $\sup_{B \in \mathcal{D}}$, it follows that

$$\|C_n^Q\| - M_n \leq \|C_n\| \leq \|C_n^Q\| + M_n, \qquad P\text{-a.s.}$$

We first assume $f$ to be bounded. Since $\|C_n^Q\| \to N$ stably under $Q$, given a bounded random variable $Z$ on $(S^\infty, \mathcal{B}^\infty)$, one obtains

$$\int \phi(\|C_n^Q\|) Z \, dQ \longrightarrow \int N(\phi) Z \, dQ$$

for each bounded continuous $\phi \colon \mathbb{R} \to \mathbb{R}$, where $N(\phi) = \int \phi(x) N(\cdot, dx)$.

Letting $Z = fI_H/P(H)$ with $H \in \mathcal{B}^\infty$ and $P(H) > 0$, it follows that $\|C_n^Q\| \to N$ stably under $P$. Therefore, it suffices to prove that $EM_n \to 0$. Given $\varepsilon > 0$, since $\|W_n\|$ is uniformly integrable under $Q$, there exists some $c > 0$ such that

$$E_Q\{\|W_n\| I_{\{\|W_n\|>c\}}\} < \frac{\varepsilon}{\sup f} \qquad \text{for all } n.$$

Since $M_n$ is $\mathcal{G}_n$-measurable,

$$EM_n = E_Q(fM_n) = E_Q(E_Q(f|\mathcal{G}_n)M_n) = E_Q(|U_n|\|W_n\|)$$
$$\leq c E_Q|U_n| + (\sup f) E_Q(\|W_n\| I_{\{\|W_n\|>c\}}) < c E_Q|U_n| + \varepsilon \qquad \text{for all } n.$$

Therefore, the martingale convergence theorem implies that

$$\limsup_n EM_n \leq c \limsup_n E_Q|U_n| + \varepsilon = \varepsilon.$$

This concludes the proof when $f$ is bounded.

Next, let $f$ be any density. Fix $k > 0$ such that $P(f \leq k) > 0$ and define $K = \{f \leq k\}$ and $P_K(\cdot) = P(\cdot|K)$. Then, $P_K$ has the bounded density $fI_K/P(K)$ with respect to $Q$. By what has already been proven, $\|C_n^{P_K}\| \to N$ stably under $P_K$, where

$$C_n^{P_K}(B) = E_{P_K}(W_n(B)|\mathcal{G}_n) = \frac{E\{I_K W_n(B)|\mathcal{G}_n\}}{E(I_K|\mathcal{G}_n)}, \qquad P_K\text{-a.s.}$$

Letting $R_n = I_K - E(I_K|\mathcal{G}_n)$, it follows that

$$E\{I_K\|C_n - C_n^{P_k}\|\} = E\left\{I_K \sup_{B\in\mathcal{D}}\left|\frac{E\{R_n W_n(B)|\mathcal{G}_n\}}{E(I_K|\mathcal{G}_n)}\right|\right\}$$

$$\leq E\left\{I_K \frac{E\{|R_n|\|W_n\||\mathcal{G}_n\}}{E(I_K|\mathcal{G}_n)}\right\} = E\{|R_n|\|W_n\|\}$$

$$\leq cE|R_n| + E\{\|W_n\|I_{\{\|W_n\|>c\}}\} \qquad \text{for all } c > 0.$$

Since $E|R_n| \to 0$ and $\|W_n\|$ is uniformly integrable under $P$, arguing as above gives that

$$E_{P_K}|\|C_n\| - \|C_n^{P_k}\|| \leq \frac{E\{I_K\|C_n - C_n^{P_k}\|\}}{P(K)} \longrightarrow 0.$$

Therefore, $\|C_n\| \to N$ stably under $P_K$. Finally, fix $H \in \mathcal{B}^\infty$, $P(H) > 0$ and a bounded continuous function $\phi: \mathbb{R} \to \mathbb{R}$. Then $P(H \cap K) = P(H \cap \{f \leq k\}) > 0$ for $k$ sufficiently large and

$$P(H)|E(\phi(\|C_n\|)|H) - E(N(\phi)|H)|$$
$$\leq 2\sup|\phi|P(f > k) + |E(\phi(\|C_n\|)|H \cap K) - E(N(\phi)|H \cap K)|.$$

Since $E(\phi(\|C_n\|)|H \cap K) \to E(N(\phi)|H \cap K)$ as $n \to \infty$ and $P(f > k) \to 0$ as $k \to \infty$, this concludes the proof. □

Next, we deal with the particular case $Q = Q_0$, where $Q_0$ is a Ferguson–Dirichlet law on $(S^\infty, \mathcal{B}^\infty)$. If $P \ll Q_0$ with a density satisfying a certain condition, the convergence rate of $\|\mu_n - a_n\|$ can be remarkably improved.

**Theorem 2.** *Suppose $\mathcal{D}$ is countably determined and $\sup_n E_{Q_0}\|W_n\|^2 < \infty$. Then, $\sqrt{n}\|C_n\| = n\|\mu_n - a_n\|$ converges a.s., provided $P \ll Q_0$ and*

$$E_{Q_0}(f^2) - E_{Q_0}\{E_{Q_0}(f|\mathcal{G}_n)^2\} = \mathrm{O}\left(\frac{1}{n}\right) \qquad \text{for some version } f \text{ of } \frac{\mathrm{d}P}{\mathrm{d}Q_0}.$$

**Proof.** Let $D_n(B) = \sqrt{n}C_n(B)$. Then, $\|D_n\|$ is $\mathcal{G}_n$-measurable (as $\mathcal{D}$ is countably determined) and

$$E(\|D_{n+1}\|\,|\mathcal{G}_n) = E\left(\sup_{B\in\mathcal{D}}\left|\sum_{i=1}^{n+1} I_B(X_i) - (n+1)E(\mu(B)|\mathcal{G}_{n+1})\right|\,\bigg|\mathcal{G}_n\right)$$

$$\geq \sup_{B\in\mathcal{D}}\left|E\left(\sum_{i=1}^{n+1} I_B(X_i)\bigg|\mathcal{G}_n\right) - (n+1)E(\mu(B)|\mathcal{G}_n)\right|$$

$$= \sup_{B\in\mathcal{D}}\left|\sum_{i=1}^{n} I_B(X_i) - nE(\mu(B)|\mathcal{G}_n)\right| = \|D_n\| \qquad \text{a.s.}$$

Since $\|D_n\|$ is a $\mathcal{G}_n$-submartingale, it suffices to prove that $\sup_n E\|D_n\| < \infty$.

Let $U_n = f - E_0(f|\mathcal{G}_n)$, where $E_0$ stands for $E_{Q_0}$. By assumption, there exist $c_1, c_2 > 0$ such that

$$E_0\|W_n\|^2 \leq c_1, \qquad nE_0U_n^2 = n\{E_0(f^2) - E_0(E_0(f|\mathcal{G}_n)^2)\} \leq c_2 \qquad \text{for all } n.$$

As noted in Section 1, since $Q_0$ is a Ferguson–Dirichlet law, there is an $\alpha > 0$ such that

$$\sqrt{n}\|C_n^{Q_0}\| = \sqrt{n}\sup_{B \in \mathcal{D}}|E_0(W_n(B)|\mathcal{G}_n)| \leq \alpha \qquad \text{for all } n.$$

Define $M_n = \frac{E_0(|U_n|\|W_n\| \mid \mathcal{G}_n)}{E_0(f|\mathcal{G}_n)}$ and recall that $\|C_n\| \leq \|C_n^{Q_0}\| + M_n$, $P$-a.s.; see the proof of Theorem 1. Then, for all $n$, one obtains

$$E\|D_n\| = \sqrt{n}E\|C_n\| \leq \sqrt{n}(E\|C_n^{Q_0}\| + EM_n) \leq \alpha + \sqrt{n}E_0(fM_n)$$
$$= \alpha + \sqrt{n}E_0(|U_n|\|W_n\|) \leq \alpha + \sqrt{n}\sqrt{E_0U_n^2 E_0\|W_n\|^2}$$
$$\leq \alpha + \sqrt{c_1 nE_0U_n^2} \leq \alpha + \sqrt{c_1 c_2}. \qquad \square$$

Finally, we clarify a point raised in Section 1.

**Remark 3.** There is a long list of (countably determined) choices of $\mathcal{D}$ such that

$$\sup_n E\|W_n\|^k \leq c(k) \qquad \text{for all } k \geq 1, \text{ if } X \text{ is i.i.d.},$$

where $c(k)$ is some universal constant; see, for example, Sections 2.14.1 and 2.14.2 of [21]. Fix one such $\mathcal{D}$, $k \geq 1$, and suppose that $S$ is Polish and $\mathcal{B}$ is the Borel $\sigma$-field. If $X$ is exchangeable, then de Finetti's theorem yields $E(\|W_n\|^k | \mathcal{T}) \leq c(k)$ a.s. for all $n$, where $\mathcal{T}$ is the tail $\sigma$-field of $X$. Hence, $E\|W_n\|^k = E\{E(\|W_n\|^k | \mathcal{T})\} \leq c(k)$ for all $n$. This proves inequality (2).

## 3. Exchangeable data with finite state space

When $X$ is exchangeable and $S$ finite, there is some overlap between Theorem 1 and a result of Bernstein and von Mises.

### 3.1. Connections with the Bernstein–von Mises theorem

For each $\theta$ in an open set $\Theta \subset \mathbb{R}^k$, let $P_\theta$ be a product probability on $(S^\infty, \mathcal{B}^\infty)$ (that is, $X$ is i.i.d. under $P_\theta$). Suppose the map $\theta \mapsto P_\theta(B)$ is Borel measurable for fixed $B \in \mathcal{B}^\infty$. Given a (prior) probability $\pi$ on the Borel subsets of $\Theta$, define

$$P(B) = \int P_\theta(B)\pi(\mathrm{d}\theta), \qquad B \in \mathcal{B}^\infty.$$

Roughly speaking, the Bernstein–von Mises (BvM) theorem can be stated as follows. Suppose $\pi$ is absolutely continuous with respect to Lebesgue measure and the statistical model $(P_\theta : \theta \in \Theta)$ is suitably "smooth" (we refer to [13] for a detailed exposition of what "smooth" means). For each $n$, suppose that $\theta$ admits a (consistent) maximum likelihood estimator $\widehat{\theta}_n$. Further, suppose the prior $\pi$ possesses the first moment and denote by $\theta_n^*$ the posterior mean of $\theta$. Then,

$$\sqrt{n}(\widehat{\theta}_n - \theta_n^*) \xrightarrow{P_{\theta_0}} 0$$

for each $\theta_0 \in \Theta$ such that the density of $\pi$ is strictly positive and continuous at $\theta_0$.

Actually, the BvM theorem yields much more than asserted; what is reported above is just the corollary connected to this paper. We refer to [13] and [14] for more information and historical notes; see also [18].

Assuming a smooth, finite-dimensional statistical model is fundamental; see, for example, [11]. Indeed, the BvM theorem does not apply when the only information is that $X$ is exchangeable (or even c.i.d.) and $P \ll Q$ for some reference probability $Q$. One exception, however, is when $S$ is finite.

Let us suppose

$$S = \{x_1, \ldots, x_k, x_{k+1}\}, \qquad X \text{ is exchangeable}, \qquad P(X_1 = x) > 0$$
$$\text{for all } x \in S \text{ and } \mathcal{D} = \mathcal{B} = \text{ power set of } S.$$

Also, let $\lambda$ denote Lebesgue measure on $\mathbb{R}^k$ and $\pi$ the probability distribution of

$$\theta = (\mu\{x_1\}, \ldots, \mu\{x_k\}).$$

As noted in Section 1, $\pi \ll \lambda$ if and only if $P \ll Q_0$ for some choice of $Q_0$. Since $\mathcal{D}$ is finite and $X$ exchangeable under $P$ and $Q_0$, $\|W_n\|$ is uniformly integrable under $P$ and $Q_0$. Thus, Theorem 1 yields $\|C_n\| \xrightarrow{P} 0$ whenever $\pi \ll \lambda$. On the other hand, $\pi$ is the prior distribution for this problem. The underlying statistical model is smooth and finite-dimensional (it is just a multinomial model). Further, for each $n$, the maximum likelihood estimator and the posterior mean of $\theta$ are, respectively,

$$\widehat{\theta}_n = (\mu_n\{x_1\}, \ldots, \mu_n\{x_k\}), \qquad \theta_n^* = (a_n\{x_1\}, \ldots, a_n\{x_k\}).$$

Thus, the BvM theorem implies that $\|C_n\| \xrightarrow{P} 0$, provided $\pi \ll \lambda$ and the density of $\pi$ is continuous on the complement of a $\pi$-null set.

To sum up, in this particular case, the same conclusions as from Theorem 1 can be drawn from the BvM theorem. Unlike the latter, however, Theorem 1 does not require any conditions on the density of $\pi$.

### 3.2. Some consequences of Theorems 1 and 2

In this subsection, we focus on $S = \{0, 1\}$. Thus, $\mathcal{D} = \mathcal{B} =$ power set of $S$ and $\lambda$ denotes Lebesgue measure on $\mathbb{R}$. Let $\mathcal{N}(0, a)$ denote the one-dimensional Gaussian law with mean 0 and variance $a \geq 0$ (where $\mathcal{N}(0,0) = \delta_0$). Our first result allows $\pi$ to have a discrete part.

**Corollary 4.** *With $S = \{0,1\}$, let $\pi$ be the probability distribution of $\mu\{1\}$ and*

$$\Delta = \{\theta \in [0,1] : \pi\{\theta\} > 0\}, \qquad A = \{\omega \in S^\infty : \mu(\omega, \{1\}) \in \Delta\}.$$

*Define the random probability measure $N$ on $\mathbb{R}$ as*

$$N = (1 - I_A)\delta_0 + I_A \mathcal{N}(0, \mu\{1\}(1 - \mu\{1\})).$$

*If $X$ is exchangeable and $\pi$ does not have a singular continuous part, then*

$$C_n\{1\} \to N \quad \text{stably} \quad \text{and} \quad \|C_n\| \to N \circ h^{-1} \quad \text{stably},$$

*where $h(x) = |x|$, $x \in \mathbb{R}$, is the modulus function.*

**Proof.** By standard arguments, the corollary holds when $\pi(\Delta) \in (0,1)$, provided it holds when $\pi(\Delta) = 0$ and $\pi(\Delta) = 1$. Let $\pi(\Delta) = 0$. Then, $\pi \ll \lambda$ as $\pi$ does not have a singular continuous part, and the corollary follows from Theorem 1. Thus, it can be assumed that $\pi(\Delta) = 1$. Since $C_n\{0\} = -C_n\{1\}$, $\|C_n\| = |C_n\{1\}|$ and the modulus function is continuous, it suffices to prove that $C_n\{1\} \to N$ stably.

Next, exchangeability of $X$ implies that $W_n\{1\} \to \mathcal{N}(0, \mu\{1\}(1 - \mu\{1\}))$ stably; see, for example, Theorem 3.1 of [6]. Since $\pi(\Delta) = 1$, we have $N = \mathcal{N}(0, \mu\{1\}(1 - \mu\{1\}))$ a.s. Hence, it is enough to show that $E|C_n\{1\} - W_n\{1\}| \to 0$.

Fix $\varepsilon > 0$ and let $M_n = W_n\{1\}$. Since $X$ is exchangeable, $M_n$ is uniformly integrable. Therefore, there exists some $c > 0$ such that

$$\sup_n E(|M_n|I_{\{|M_n|>c\}}) < \frac{\varepsilon}{4}.$$

Define $\phi(x) = x$ if $|x| \leq c$, $\phi(x) = c$ if $x > c$, and $\phi(x) = -c$ if $x < -c$. Since $C_n\{1\} = E(M_n|\mathcal{G}_n)$ a.s., it follows that

$$\begin{aligned}
E|C_n\{1\} - W_n\{1\}| &\leq E|E(M_n|\mathcal{G}_n) - E(\phi(M_n)|\mathcal{G}_n)| \\
&\quad + E|E(\phi(M_n)|\mathcal{G}_n) - \phi(M_n)| + E|\phi(M_n) - M_n| \\
&\leq E|E(\phi(M_n)|\mathcal{G}_n) - \phi(M_n)| + 4E(|M_n|I_{\{|M_n|>c\}}) \\
&< E|E(\phi(M_n)|\mathcal{G}_n) - \phi(M_n)| + \varepsilon \qquad \text{for all } n.
\end{aligned}$$

Write $\Delta = \{a_1, a_2, \ldots\}$ and $M_{n,j} = \sqrt{n}(\mu_n\{1\} - a_j)$. Since $\sigma(M_{n,j}) \subset \mathcal{G}_n$ and $P(\mu\{1\} \in \Delta) = \pi(\Delta) = 1$, one also obtains

$$\begin{aligned}
&E|E(\phi(M_n)|\mathcal{G}_n) - \phi(M_n)| \\
&= \sum_j E|E(\phi(M_{n,j})I_{\{\mu\{1\}=a_j\}}|\mathcal{G}_n) - \phi(M_{n,j})I_{\{\mu\{1\}=a_j\}}|
\end{aligned}$$

$$= \sum_j E|\phi(M_{n,j})\{P(\mu\{1\} = a_j|\mathcal{G}_n) - I_{\{\mu\{1\}=a_j\}}\}|$$

$$\le c \sum_{j=1}^m E|P(\mu\{1\} = a_j|\mathcal{G}_n) - I_{\{\mu\{1\}=a_j\}}| + 2c \sum_{j>m} \pi\{a_j\} \qquad \text{for all } m, n.$$

By the martingale convergence theorem, $E|P(\mu\{1\} = a_j|\mathcal{G}_n) - I_{\{\mu\{1\}=a_j\}}| \to 0$ as $n \to \infty$, for each $j$. Thus,

$$\limsup_n E|C_n\{1\} - W_n\{1\}| \le \varepsilon + 2c \sum_{j>m} \pi\{a_j\} \qquad \text{for all } m.$$

Taking the limit as $m \to \infty$ completes the proof. $\square$

If $\pi$ is singular continuous, we conjecture that $C_n\{1\}$ converges stably to a non-null limit. However, we do not have a proof.

In the next result, a real function $g$ on $(0,1)$ is said to be *almost Lipschitz* in the case where $x \mapsto g(x)x^a(1-x)^b$ is Lipschitz on $(0,1)$ for some reals $a, b < 1$.

**Corollary 5.** *Suppose $S = \{0, 1\}$, $X$ is exchangeable and $\pi$ is the probability distribution of $\mu\{1\}$. If $\pi$ admits an almost Lipschitz density with respect to $\lambda$, then $\sqrt{n}\|C_n\|$ converges a.s. to a real random variable.*

**Proof.** Let $V = \mu\{1\}$. By assumption, there exist $a, b < 1$ and a version $g$ of $\frac{d\pi}{d\lambda}$ such that $\phi(\theta) = g(\theta)\theta^a(1-\theta)^b$ is Lipschitz on $(0,1)$. For each $u_1, u_2 > 0$, we can take $Q_0$ such that $V$ has a beta-distribution with parameters $u_1, u_2$ under $Q_0$. Let $Q_0$ be such that $V$ has a beta-distribution with parameters $u_1 = 1 - a$ and $u_2 = 1 - b$ under $Q_0$. Then, for any $n \ge 1$ and $x_1, \ldots, x_n \in \{0, 1\}$, one obtains

$$P(X_1 = x_1, \ldots, X_n = x_n)$$
$$= \int_0^1 \theta^r (1-\theta)^{n-r} \pi(d\theta)$$
$$= \int_0^1 \theta^{r-a}(1-\theta)^{n-r-b}\phi(\theta)\,d\theta$$
$$= c \int V^r(1-V)^{n-r}\phi(V)\,dQ_0, \qquad \text{where } r = \sum_{i=1}^n x_i \text{ and } c > 0 \text{ is a constant.}$$

Let $h = c\phi$. Then, $h$ is Lipschitz and $f = h(V)$ is a version of $\frac{dP}{dQ_0}$.

Let $V_n = E_0(V|\mathcal{G}_n)$, where $E_0$ stands for $E_{Q_0}$. Since $h$ is Lipschitz,

$$|f - E_0(f|\mathcal{G}_n)| \le |h(V) - h(V_n)| + E_0(|h(V) - h(V_n)||\mathcal{G}_n)$$
$$\le d|V - V_n| + dE_0(|V - V_n||\mathcal{G}_n),$$

where $d$ is the Lipschitz constant of $h$. Since $E_0\|C_n^{Q_0}\|^2 \leq E_0\|W_n\|^2$ and

$$\sqrt{n}|V - V_n| = |C_n^{Q_0}\{1\} - W_n\{1\}| \leq \|C_n^{Q_0}\| + \|W_n\|,$$

it follows that

$$E_0(f^2) - E_0(E_0(f|\mathcal{G}_n)^2) = E_0\{(f - E_0(f|\mathcal{G}_n))^2\} \leq 4d^2 E_0\{(V - V_n)^2\}$$

$$\leq \frac{4d^2}{n} E_0\{(\|C_n^{Q_0}\| + \|W_n\|)^2\} \leq \frac{16d^2}{n} E_0\|W_n\|^2.$$

Since $\sup_n E_0\|W_n\|^2 < \infty$, we have $E_0(f^2) - E_0(E_0(f|\mathcal{G}_n)^2) = \mathrm{O}(1/n)$. An application of Theorem 2 completes the proof. □

Corollaries 4 and 5 deal with $S = \{0, 1\}$, but similar results can be proven for any finite $S$; see also [12] and [19].

## 4. Generalized Pólya urns

In this section, based on Examples 1.3 and 3.5 of [6], the asymptotic behavior of $\|C_n\|$ is investigated for a certain c.i.d. sequence.

Let $(\mathcal{Y}, \mathcal{B}_\mathcal{Y})$ be a measurable space, $\mathcal{B}_+$ the Borel $\sigma$-field on $(0, \infty)$ and

$$S = \mathcal{Y} \times (0, \infty), \qquad \mathcal{B} = \mathcal{B}_\mathcal{Y} \otimes \mathcal{B}_+, \qquad X_n = (Y_n, Z_n),$$

where $Y_n(\omega) = y_n, Z_n(\omega) = z_n$ for all $\omega = (y_1, z_1, y_2, z_2, \ldots) \in S^\infty$.

Given a law $P$ on $\mathcal{B}^\infty$, it is assumed that

$$P(Y_{n+1} \in B|\mathcal{G}_n) = \frac{\alpha P(Y_1 \in B) + \sum_{i=1}^n Z_i I_B(Y_i)}{\alpha + \sum_{i=1}^n Z_i} \qquad \text{a.s., } n \geq 1, \tag{4}$$

$$P(Z_{n+1} \in C|X_1, \ldots, X_n, Y_{n+1}) = P(Z_1 \in C) \qquad \text{a.s., } n \geq 0, \tag{5}$$

for some constant $\alpha > 0$ and all $B \in \mathcal{B}_\mathcal{Y}, C \in \mathcal{B}_+$. Note that $(Z_n)$ is i.i.d. and $Z_{n+1}$ is independent of $(Y_1, Z_1, \ldots, Y_n, Z_n, Y_{n+1})$ for all $n \geq 0$.

In real problems, the $Z_n$ should be viewed as weights, while the $Y_n$ describe the phenomenon of interest. As an example, consider an urn containing white and black balls. At each time $n \geq 1$, a ball is drawn and then replaced together with $Z_n$ more balls of the same color. Let $Y_n$ be the indicator of the event {white ball at time $n$} and suppose that $Z_n$ is chosen according to a fixed distribution on the integers, independently of $(Y_1, Z_1, \ldots, Y_{n-1}, Z_{n-1}, Y_n)$. The predictive distributions of $X$ are then given by (4)–(5). Also, note that the probability law of $(Y_n)$ is Ferguson–Dirichlet in the case where $Z_n = 1$ for all $n$.

It is not hard to prove that $X$ is c.i.d. We state this fact as a lemma.

**Lemma 6.** *The sequence $X$ assessed according to (4)–(5) is c.i.d.*

**Proof.** Fix $k > n \geq 0$ and $A \in \mathcal{B}_\mathcal{Y} \otimes \mathcal{B}_+$. By a monotone class argument, it can be assumed that $A = B \times C$, where $B \in \mathcal{B}_\mathcal{Y}$ and $C \in \mathcal{B}_+$. Further, it can be assumed that $k = n + 2$. Let $n = 0$ and $\mathcal{G}_0$ be the trivial $\sigma$-field. Since $X_2 \sim X_1$ (as is easily seen), $E(I_B(Y_2)I_C(Z_2)|\mathcal{G}_0) = E(I_B(Y_1)I_C(Z_1)|\mathcal{G}_0)$ a.s. If $n \geq 1$, define

$$\mathcal{G}_n^* = \sigma(X_1, \ldots, X_n, Z_{n+1}).$$

Noting that $E(I_B(Y_{n+1})|\mathcal{G}_n^*) = E(I_B(Y_{n+1})|\mathcal{G}_n)$ a.s., one obtains

$$E(I_B(Y_{n+2})|\mathcal{G}_n^*) = E\{E(I_B(Y_{n+2})|\mathcal{G}_{n+1})|\mathcal{G}_n^*\}$$
$$= \frac{\alpha P(Y_1 \in B) + \sum_{i=1}^n Z_i I_B(Y_i) + Z_{n+1} E(I_B(Y_{n+1})|\mathcal{G}_n^*)}{\alpha + \sum_{i=1}^{n+1} Z_i}$$
$$= \frac{(\alpha + \sum_{i=1}^n Z_i) E(I_B(Y_{n+1})|\mathcal{G}_n) + Z_{n+1} E(I_B(Y_{n+1})|\mathcal{G}_n)}{\alpha + \sum_{i=1}^{n+1} Z_i}.$$
$$= E(I_B(Y_{n+1})|\mathcal{G}_n) = E(I_B(Y_{n+1})|\mathcal{G}_n^*) \quad \text{a.s.}$$

Finally, since $\mathcal{G}_n \subset \mathcal{G}_n^*$, the previous equality implies that

$$E(I_B(Y_{n+2})I_C(Z_{n+2})|\mathcal{G}_n) = P(Z_1 \in C)E\{E(I_B(Y_{n+2})|\mathcal{G}_n^*)|\mathcal{G}_n\}$$
$$= P(Z_1 \in C)E\{E(I_B(Y_{n+1})|\mathcal{G}_n^*)|\mathcal{G}_n\}$$
$$= E(I_B(Y_{n+1})I_C(Z_{n+1})|\mathcal{G}_n) \quad \text{a.s.}$$

Therefore, $X$ is c.i.d. □

Usually, one is interested in predicting $Y_n$ more than $Z_n$. Thus, in the sequel, we focus on $P(Y_{n+1} \in B|\mathcal{G}_n)$. For each $B \in \mathcal{B}_\mathcal{Y}$, we write

$$C_n(B) = C_n(B \times (0, \infty)), \qquad a_n(B) = a_n(B \times (0, \infty)) = P(Y_{n+1} \in B|\mathcal{G}_n),$$

and so on.

In Example 3.5 of [6], assuming $EZ_1^2 < \infty$, it is shown that

$$C_n(B) \to \mathcal{N}(0, \sigma_B^2) \quad \text{stably, where } \sigma_B^2 = \frac{\text{var}(Z_1)}{(EZ_1)^2}\mu(B)(1 - \mu(B)).$$

Here, we prove that $C_n$ converges stably when regarded as a map $C_n : S^\infty \to l^\infty(\mathcal{D})$, where $l^\infty(\mathcal{D})$ is the space of real bounded functions on $\mathcal{D}$ equipped with uniform distance; see Section 1.5 of [21]. In particular, stable convergence of $C_n$ as a random element of $l^\infty(\mathcal{D})$ implies stable convergence of $\|C_n\| = \sup_{B \in \mathcal{D}} |C_n(B)|$.

Intuitively, the stable limit of $C_n$ (when it exists) is connected to the Brownian bridge. Let $B_1, B_2, \ldots$ be pairwise disjoint elements of $\mathcal{B}_\mathcal{Y}$ and

$$\mathcal{D} = \{B_k \times (0, \infty) : k \geq 1\}, \qquad T_0 = 0, \qquad T_k = \sum_{i=1}^k \mu(B_i).$$

Also, let $G$ be a standard Brownian bridge process on some probability space $(\Omega_0, \mathcal{A}_0, P_0)$. For fixed $\omega \in S^\infty$,

$$L(\omega, B_k) = \frac{\sqrt{\operatorname{var}(Z_1)}}{EZ_1}\{G(T_k(\omega)) - G(T_{k-1}(\omega))\}$$

is a real random variable on $(\Omega_0, \mathcal{A}_0, P_0)$. Since the $B_k$ are pairwise disjoint and $G$ has continuous paths, $L(\omega, B_k) \to 0$ as $k \to \infty$. It thus makes sense to define $M(\omega, \cdot)$ as the probability distribution of $L(\omega) = (L(\omega, B_1), L(\omega, B_2), \ldots)$, that is,

$$M(\omega, A) = P_0(L(\omega) \in A) \qquad \text{for each Borel set } A \subset l^\infty(\mathcal{D}).$$

Similarly, let $N(\omega, \cdot)$ be the probability distribution of $\sup_{k \geq 1} |L(\omega, B_k)|$, that is,

$$N(\omega, A) = P_0\left(\sup_{k \geq 1} |L(\omega, B_k)| \in A\right) \qquad \text{for each Borel set } A \subset \mathbb{R}.$$

**Theorem 7.** *Suppose $B_1, B_2, \ldots \in \mathcal{B}_\mathcal{Y}$ are pairwise disjoint and $\mathcal{D}$, $M$, $N$ are defined as above. Let $X$ be assessed according to (4)–(5) with $a \leq Z_1 \leq b$ a.s. for some constants $0 < a < b$. Then,*

$$\sup_n E\|W_n\|^2 \leq c\sqrt{P\left(Y_1 \in \bigcup_k B_k\right)} \tag{6}$$

*for some constant $c$ independent of the $B_k$, and $C_n \to M$ stably (in the metric space $l^\infty(\mathcal{D})$). In particular, $\|C_n\| \to N$ stably.*

Let $Q_1$ denote the probability law of a sequence $X$ satisfying (4)–(5) and $a \leq Z_1 \leq b$ a.s. In view of Theorem 7, $Q_1$ can play the role of $Q$ in Theorem 1. That is, for an arbitrary c.i.d. sequence $X$ with distribution $P$, one has $\|C_n\| \to N$ stably, provided $P \ll Q_1$ and $\|W_n\|$ is uniformly integrable under $P$. The condition of pairwise disjoint $B_k$ is actually rather strong. However, it holds in at least two relevant situations: when a single set $B$ is involved, and when $S = \{x_1, x_2, \ldots\}$ is countable and $B_k = \{x_k\}$ for all $k$.

**Proof of Theorem 7.** This proof involves some simple but long calculations. Accordingly, we provide only a sketch of the proof and refer to [7] for details.

Since $X$ is c.i.d., for fixed $B \in \mathcal{B}_\mathcal{Y}$, one has $a_n(B) = E(\mu(B)|\mathcal{G}_n)$ a.s. Hence, $(a_n(B) : n \geq 1)$ is a $\mathcal{G}_n$-martingale with $a_n(B) \xrightarrow{\text{a.s.}} \mu(B)$ and this implies that

$$E\{(a_{n+1}(B) - \mu(B))^2\} = E\left\{\left(\sum_{j > n}(a_j(B) - a_{j+1}(B))\right)^2\right\} = \sum_{j > n} E\{(a_j(B) - a_{j+1}(B))^2\}.$$

Replacing $a_j(B)$ by (4) and using the fact that $a \leq Z_i \leq b$ a.s. for all $i$, a long but straightforward calculation yields $\sum_{j > n} E\{(a_j(B) - a_{j+1}(B))^2\} \leq \frac{c_1}{n} P(Y_1 \in B)$, where

$c_1$ is a constant independent of $B$. It follows that

$$E\|a_{n+1} - \mu\|^2 = E\Big\{\sup_k(a_{n+1}(B_k) - \mu(B_k))^2\Big\} \le \sum_k E\{(a_{n+1}(B_k) - \mu(B_k))^2\}$$

$$= \sum_k \sum_{j>n} E\{(a_j(B_k) - a_{j+1}(B_k))^2\} \le \frac{c_1}{n}\sum_k P(Y_1 \in B_k)$$

$$= \frac{c_1}{n} P\Big(Y_1 \in \bigcup_k B_k\Big) \qquad \text{as the } B_k \text{ are pairwise disjoint.}$$

Precisely as above, after some algebra, one obtains

$$E\|\mu_n - a_{n+1}\|^2 \le \frac{c_2}{n}\sqrt{P\Big(Y_1 \in \bigcup_k B_k\Big)}$$

for some constant $c_2$ independent of $B_1, B_2, \ldots$. Therefore,

$$E\|W_n\|^2 = nE\|\mu_n - \mu\|^2 \le 2nE\|\mu_n - a_{n+1}\|^2 + 2nE\|a_{n+1} - \mu\|^2 \le c\sqrt{P\Big(Y_1 \in \bigcup_k B_k\Big)},$$

where $c = 2(c_1 + c_2)$. This proves inequality (6).

It remains to prove that $C_n \to M$ stably (in the metric space $l^\infty(\mathcal{D})$). For each $m \ge 1$, let $\Sigma_m$ be the $m \times m$ matrix with elements

$$\sigma_{k,j} = \frac{\operatorname{var}(Z_1)}{(EZ_1)^2}(\mu(B_k \cap B_j) - \mu(B_k)\mu(B_j)), \qquad k, j = 1, \ldots, m.$$

By Theorems 1.5.4 and 1.5.6 of [21], for $C_n \to M$ stably, it is enough that:

(i) (finite-dimensional convergence):

$$(C_n(B_1), \ldots, C_n(B_m)) \to \mathcal{N}_m(0, \Sigma_m) \qquad \text{stably for each } m \ge 1,$$

where $\mathcal{N}_m(0, \Sigma_m)$ is the $m$-dimensional Gaussian law with mean 0 and covariance matrix $\Sigma_m$;

(ii) (asymptotic tightness): for each $\varepsilon, \delta > 0$, there exists some $m \ge 1$ such that

$$\limsup_n P\Big(\sup_{r,s>m}|C_n(B_r) - C_n(B_s)| > \varepsilon\Big) < \delta.$$

Fix $m \ge 1$, $b_1, \ldots, b_m \in \mathbb{R}$ and define $R_n = \sum_{k=1}^m b_k I_{B_k}(Y_n)$. Since $(R_n : n \ge 1)$ is c.i.d., arguing exactly as in Example 3.5 of [6], one obtains

$$\sum_{k=1}^m b_k C_n(B_k) = \frac{\sum_{i=1}^n \{R_i - E(R_{n+1}|\mathcal{G}_n)\}}{\sqrt{n}} \longrightarrow \mathcal{N}\Big(0, \sum_{k,j} b_k b_j \sigma_{k,j}\Big) \qquad \text{stably.}$$

Since $b_1, \ldots, b_m$ are arbitrary, (i) holds. To check (ii), given $\varepsilon, \delta > 0$, take $m$ such that

$$P\left(Y_1 \in \bigcup_{r>m} B_r\right) < \left(\frac{\varepsilon^2 \delta}{4c}\right)^2,$$

where $c$ is the constant involved in (6). By what has already been proven,

$$P\left(\sup_{r,s>m} |C_n(B_r) - C_n(B_s)| > \varepsilon\right) \leq P\left(2 \sup_{r>m} |C_n(B_r)| > \varepsilon\right)$$
$$\leq P\left(2E\left(\sup_{r>m} |W_n(B_r)| | \mathcal{G}_n\right) > \varepsilon\right) \leq \frac{4}{\varepsilon^2} E\left\{\sup_{r>m} W_n(B_r)^2\right\}$$
$$\leq \frac{4c}{\varepsilon^2} \sqrt{P\left(Y_1 \in \bigcup_{r>m} B_r\right)} < \delta.$$

Thus, (ii) holds and this completes the proof. $\square$

## Acknowledgments

This paper benefited from the helpful suggestions of two anonymous referees.